\documentclass[A4,12pt]{amsart}            
\usepackage{amsmath,amssymb,amsrefs,amsfonts,eucal,yfonts,amscd,comment}

\newtheorem{theorem}{Theorem}[section]
\newtheorem{proposition}[theorem]{Proposition}
\newtheorem{lemma}[theorem]{Lemma}
\newtheorem{corollary}[theorem]{Corollary}

\theoremstyle{definition}

\newtheorem{remark}[theorem]{Remark}
\newtheorem{remarks}[theorem]{Remarks}

\newcommand{\Real}{\mathbb R}                                   

\def\ftp{\ensuremath{\overline{\otimes}}}

\begin{document}
\title{Tensor products of $f$-algebras}

\author{G.J.H.M. Buskes}
\address{Department of Mathematics, University of Mississippi, University, MS 38677, USA}\email{mmbuskes@olemiss.edu}
\author{A.W. Wickstead}
\address{Pure Mathematics Research Centre, Queen's University
Belfast, Belfast BT7 1NN, Northern Ireland, UK} \email{A.Wickstead@qub.ac.uk}

\begin{abstract}
We construct the tensor product for $f$-algebras, including proving a universal property for it, and investigate how it preserves algebraic properties of the factors.
\end{abstract}

\thanks{The authors would like to thank the referee for several helpful suggestions and, in particular, for bringing \cite{T} to their attention.}
 \maketitle
\section{Introduction}
An $f$-algebra is a vector lattice $A$ with an associative  multiplication $\star$ such that
\begin{enumerate}
\item $a,b\in A_+\Rightarrow a\star b\in A_+$
\item If $a,b,c\in A_+$ with $a\wedge b=0$ then $(a\star c)\wedge b=(c\star a)\wedge b=0$.
\end{enumerate}
We are only concerned in this paper with Archimedean $f$-algebras so \emph{henceforth all $f$-algebras, and indeed all vector lattices, are assumed to be Archimedean.}

In this paper, we construct the tensor product of two $f$-algebras $A$ and $B$. To be precise we show that there is a natural way in which the Fremlin Archimedean vector lattice tensor product $A\ftp B$ can be endowed with an $f$-algebra multiplication. This tensor product has the expected universal property, at least as far as order bounded algebra homomorphisms are concerned. We also show that, unlike the case for order theoretic properties, algebraic properties behave well under this product. Again, let us emphasise here that \emph{algebra homomorphisms are not assumed to map an identity to an identity.}

Our approach is very much representational, however much that might annoy some purists. In \S 2 we derive the representations needed, which are very mild improvements on results already in the literature.  In the proof of our main result a first attempt at a proof ran into the problem that if $f\in C^\infty(\Sigma)$ and $g\in C^\infty(\Omega)$ (continuous functions into $\Real\cup\{-\infty,\infty\}$ that are finite on a dense subset) then we do not have $(\sigma,\omega)\mapsto f(\sigma)g(\omega)$ lying in $C^\infty(\Sigma\times \Omega)$. We get round this problem by working with representations in the space of continuous real-valued functions defined on dense open subsets, modulo equality on intersection of domains.

Although we can find no direct connection, we feel that we should pay homage to the first use of $f$-algebras in connection with tensor products in \cite{GL}. There is overlap between our results and those obtained (for unital and semi-prime $f$-algebras, and using different methods) in \cite{AAJ}.

\section{Representations of $f$-algebras}\label{reps}

There are many representation results in the literature for $f$-algebras, e.g. \cites{BP,BvR,HeJ,J}. We are able to prove a result that is slightly more general than any of these, but our main reason for approaching the topic of representations is that we actually need to represent the multiplication in situations when we have a functional representation of only some fragment of the $f$-algebra.

If $\Sigma$ is a topological space, then $C^\infty(\Sigma)$ is the set of continuous functions into the usual two point compactification of the reals which are finite on a dense subset of $\Sigma$. This is a lattice for the pointwise order but is not in general a vector space for pointwise addition and scalar multiplication. There may, however, be many vector subspaces. The proofs of Lemma 2.4 and Theorem 2.5 of \cite{W} suffice to prove the following result, so we do not repeat them here.

\begin{proposition}\label{orthorep}Let $\Sigma$ be a topological space, $E\subset C^\infty(\Sigma)$ a vector lattice and $T:E\to C^\infty(\Sigma)$ be such that:
\begin{enumerate}
\item  If $x\ge 0$ then $Tx\ge 0$.
\item For all $\alpha,\beta\in\Real$ and $x,y\in E$, $\alpha Tx+\beta Ty$ is defined in $C^\infty(\Sigma)$ and is equal to $T(\alpha x+\beta y)$.
\item If $x,y\in E$ with $x\perp y$ implies that $Tx\perp y$.
\end{enumerate}
Set $\Upsilon=\{\sigma\in \Sigma:\exists x\in E\text{ with } 0<|x(\sigma)|<\infty\}$ then there is $q\in C^\infty(\Upsilon)$ such that for any $x\in E$
\[Tx(\upsilon)=q(\upsilon) x(\upsilon)\]
for any $\upsilon\in \Upsilon$ for which the product is defined.
\end{proposition}

\begin{remarks}
Compared with Theorem 2.5 of \cite{W}, the differences in this result are that the operator $T$ takes values only in $C^\infty(\Sigma)$ rather than in $E$ and that the set $\Upsilon$ is not assumed to be the whole of $\Sigma$. As a consequence we have no information about $Tx$ on $\Sigma\setminus \Upsilon$. The price we pay for this generalization is that we must assume that $T$ is positive rather than order bounded.
\end{remarks}

\begin{proposition}\label{fpartrep}
 Let $(A,\star)$ be an $f$-algebra, $H$ a vector sublattice  of $A$ and $G$ a vector sublattice of $H$ such that $G\star G\subseteq H$. Let $x\mapsto \hat{x}$ be a representation of $H$ in some $C^\infty(\Sigma)$ and let $\Upsilon=\{\sigma\in\Sigma:\exists x\in G\text{ with }0<|\hat{x}(\sigma)|<\infty\}$. Then there is $w\in C^\infty(\Upsilon)_+$ such that, for all $x,y\in G$,
\[\widehat{(x\star y)}(\upsilon)=w(\upsilon) \hat{x}(\upsilon) \hat{y}(\upsilon)\]
at all points $\upsilon\in\Upsilon$ for which  the pointwise product is defined.

\begin{proof}
We will identify $H$ with a sublattice of $C^\infty(\Sigma)$ and suppress the $x\mapsto \hat{x}$ notation during this proof.
 For each $x\in G_+$ the map $M_x:y\mapsto x\star y:G\to H\subseteq C^\infty(\Sigma)$  satisfies the hypotheses of Proposition \ref{orthorep}, so there is $q_x\in C^\infty(\Upsilon)$ with $(x\star y)(\upsilon)=(M_x y)(\upsilon)=q_x(\upsilon) y(\upsilon)$ for all $\upsilon\in\Upsilon$ for which  the pointwise product is defined. For general $x\in G$ this argument can be applied to $x^+$ and $x^-$ separately to obtain $q_{x^+}$ and $q_{x^-}$ in $C^\infty(\Upsilon)$. As $x^+\perp x^-$, $M_{x^+}y\perp M_{x^-}y$ for all $y\in G$. Picking any $\upsilon\in\Upsilon$ at which both $q_{x^+}$ and $q_{x^-}$ are finite, we may find $y\in G_+$ with $0<y(\sigma)<\infty$ and we see that at least one of $q_{x^+}$ and $q_{x^-}$ must vanish at $\upsilon$. By continuity, $q_{x^+}\perp q_{x^-}$ and we may form $q_x=q_{x^+}-q_{x^-}\in C^\infty(\Upsilon)$. It is routine to verify that we now have the equality, $(x\star y)(\upsilon)=q_x(\upsilon) y(\upsilon)$ for all $\upsilon\in\Upsilon$ for which  the pointwise product is defined, for all $x\in G$.

 Now consider the map $Q:x\mapsto q_x$ mapping $G$ into $C^\infty(\Upsilon)$. From this point on, we regard $G$ as a subspace of $C^\infty(\Upsilon)$ in the obvious manner. We show that $Q$ is both linear and positive. It is clear that $Q$ respects multiplication by reals. If $x,x'\in G$, in order to prove that $q_{x+x'}=q_x+q_{x'}$ it suffices to prove that this equality holds on the dense subset of $\Upsilon$ consisting of points $\upsilon$ for which all three functions are finite. But if we choose $y\in G_+$ with $0<y(\upsilon)<\infty$ then the equality
 \begin{align*}
   q_{x+x'}(\upsilon) y(\upsilon)&=\big((x+x')\star y\big)(\upsilon)=(x\star y)(\upsilon)+(x'\star y)(\upsilon)\\
   &=q_x(\upsilon) y(\upsilon)+q_{x'}(\upsilon) y(\upsilon)=\big(q_x(\upsilon)+q_{x'}(\upsilon)\big)y(\upsilon)
 \end{align*}
makes that clear. A similar argument shows that $Q$ is positive.
 If $\upsilon\in\Upsilon$, pick $y\in G_+$ with $0<y(\upsilon)<\infty$.

We next show that if $x,y\in G$ and $x\perp y$ then $Qx\perp y$. Suppose that $y(\upsilon)\ne 0$ and  use the admissibility of $G$ to find $z\in A$ with $0<z(\upsilon)<\infty$. There is a neighbourhood $U$ of $\upsilon$ on which both sets of inequalities persist. As $x \perp y$, $x\star z\perp y$. For $\tau$ in the  the dense subset of $\Upsilon$ where the pointwise product is defined, $Qx(\tau) z(\tau)=(x\star z)(\tau)$. As $x\star z\perp y$, $Qx(\tau) z(\tau)=0$ on a dense subset of $U$ and hence $Qx(\tau)=0$ on that dense subset of $U$. By continuity, $Qx\equiv 0$ on $U$ and in particular $Qx(\upsilon)=0$.

We may thus apply Proposition \ref{orthorep} again to $Q$ to obtain $w\in C^\infty(\Upsilon)_+$ such that $Qx(\upsilon)=q_x(\upsilon)=w(\upsilon) x(\upsilon)$ whenever the product is defined. Hence, for all $x,y\in G$ we have
\[(x\star y)(\upsilon)=q_x(\upsilon) y(\upsilon)=w(\upsilon) x(\upsilon) y(\upsilon)\]
for all $\upsilon\in\Upsilon$ for which the pointwise product is defined.
\end{proof}
\end{proposition}

If the sublattice $H$ has an order unit then we have a rather simpler representation. The representation mentioned here exists by the Kakutani representation.

\begin{corollary}\label{fpartrep2}
 Let $(A,\star)$ be an $f$-algebra, $H$ a sublattice of $A$ with an order unit and $G$ a sublattice of $H$ such that $G\star G\subseteq H$. Let $x\mapsto \hat{x}$ be a representation of $H$ in some $C(\Sigma)$, where $\Sigma$ is compact Hausdorff, and let $\Upsilon=\{\sigma\in\Sigma:\exists x\in G\text{ with }0<|\hat{x}(\sigma)|\}$. Then there is $w\in C(\Upsilon)_+$ such that, for all $x,y\in G$ and all $\upsilon\in\Upsilon$,
\[\widehat{(x\star y)}(\upsilon)=w(\upsilon) \hat{x}(\upsilon) \hat{y}(\upsilon).\]
\begin{proof}
We need only note that $w$ is locally a quotient of $Tx$ by the corresponding (non-vanishing) $x$ to see that $w$ is real-valued.
\end{proof}
\end{corollary}

In particular we have the following consequence of Proposition \ref{fpartrep} giving a global representation of $f$-algebras.

\begin{corollary}\label{frep}
Let $\Sigma$ be a topological space and $A\subset C^\infty(\Sigma)$ is a vector lattice such that for all $\sigma\in\Sigma$ there is $x\in A$ with $0<|x(\sigma)|<\infty$. If $\star$ is an $f$-algebra multiplication on $A$ then there is $w\in C^\infty(\Sigma)_+$ such that, for all $x,y\in A$,
\[(x\star y)(\sigma)=w(\sigma) x(\sigma) y(\sigma)\]
for all $\sigma\in \Sigma$ for which the pointwise product is defined.

Conversely, if there is $w\in C^\infty(\Sigma)$ such that for all $x,y\in A$ there is an element $x\star y$ of $A$ with $(x\star y)(\sigma)=w(\sigma) x(\sigma) y(\sigma)$ whenever that pointwise product is defined, then $\star$ is an $f$-algebra multiplication on $A$.
\begin{proof}Use Proposition \ref{fpartrep} taking $G=H=A$ and note that in this case $\Upsilon=\Sigma$. The converse is elementary.
\end{proof}
\end{corollary}

Of course, every Archimedean vector lattice possesses many representations as such sublattices of $C^\infty(\Sigma)$. Theorem 6 of \cite{BvR} proves this corollary under the additional assumption that $\Sigma$ is Stonean. Specialising somewhat, we have:

\begin{corollary}\label{semiprimerep}
For every semi-prime $f$-algebra $(A,\star)$ there is a vector lattice representation of $A$ in some $C^\infty(\Upsilon)$, with image a vector lattice, and a strictly positive $w\in C^\infty(\Upsilon)_+$ such that, for all $x,y\in A$,
\[(x\star y)(\upsilon)=w(\upsilon) x(\upsilon) y(\upsilon)\]
for all $\upsilon\in \Upsilon$ for which the pointwise product is defined.
\begin{proof}
Start with any suitable representation of $A$ in $C^\infty(\Sigma)$ and produce $w$ as in Corollary \ref{frep}, then set $\Upsilon=\{\sigma\in\Sigma:w(\sigma)>0\}$. The restriction of the representation on $\Sigma$ to $\Upsilon$ will do what we want, as soon as we show that this restriction is one-to-one. But if $x\in A$ and $x_{|\Upsilon}=0$ then $w(\sigma) x(\sigma)^2=0$ for all $\sigma\in \Sigma$, which says that $x\star x=0$ and $x=0$ as $(A,\star)$ is semi-prime.
\end{proof}
\end{corollary}

For the work below, Corollary \ref{frep} is not quite good enough. We will work with representations not in a space $C^\infty(\Sigma)$, which as we have already commented is not a vector space in any natural way, but in the space that we will denote by $S(\Sigma)$. This is the space of continuous real valued functions defined on dense open subsets of $\Sigma$, modulo the equivalence relation that $f\sim g$ if $f$ and $g$ coincide on the intersection of their domains. There is a natural embedding $\pi$ of $C^\infty(\Sigma)$ into $S(\Sigma)$ obtained by restricting $f\in C^\infty(\Sigma)$ to the subset of $\Sigma$ on which it is finite. However, the set $S(\Sigma)$ is strictly larger and in particular is easily verified to be a vector space under the pointwise operations on the intersections of domains. In particular, $S(\Sigma)$ is a vector lattice under the pointwise order. By taking a representation of an Archimedean vector lattice $E$ in  $C^\infty(\Sigma)$ and composing with the embedding $\pi$ we obtain a representation of $E$ in $S(\Sigma)$. For such representations, Corollary \ref{frep} takes the following form.

\begin{corollary}\label{frep2}
 Let $A$ be an admissible vector sublattice of $S(\Sigma)$ for some topological space $\Sigma$. If $\star$ is an $f$-algebra multiplication on $A$ then there is $w\in S(\Sigma)_+$ such that, for all $x,y\in A$,
\[(x\star y)(\sigma)=w(\sigma) x(\sigma) y(\sigma)\eqno{(\dag)}\]
on the intersection of all the relevant domains.

Conversely, if there is $w\in S(\Sigma)_+$ such that for all $x,y\in A$ there is an element $x\star y$ of $A$ satisfying {\rm($\dag$)} on the intersection of the domains, then $\star$ is an $f$-algebra multiplication on $A$.

In particular, for any $w\in S(\Sigma)_+$, we may define an $f$-algebra multiplication on the whole of $S(\Sigma)$ by {\rm($\dag$)} on the intersection of the domains of $w,x$ and $y$.
\end{corollary}

It is immediate that if the $w$ in Corollary \ref{frep2} is strictly positive then $(A,\star)$ is semi-prime. Conversely, we have as a consequence of Corollary \ref{semiprimerep}:

\begin{corollary}\label{semiprimerep2}
For every semi-prime $f$-algebra $(A,\star)$ there is a vector lattice representation of $A$ as a vector sublattice of some $S(\Upsilon)$ and a strictly positive $w\in S(\Upsilon)_+$ such that, for all $x,y\in A$,
\[(x\star y)(\upsilon)=w(\upsilon) x(\upsilon) y(\upsilon)\]
on the intersection of all the relevant domains.
\begin{proof}
Start with any suitable representation of $A$ in $C^\infty(\Sigma)$ and produce $w$ as in Corollary \ref{frep}, then set $\Upsilon=\{\sigma\in\Sigma:w(\sigma)>0\}$. The restriction of the representation on $\Sigma$ to $\Upsilon$ will do what we want, as soon as we show that this restriction is one-to-one. But if $x\in A$ and $x_{|\Upsilon}=0$ then $w(\sigma) x(\sigma)^2=0$ for all $\sigma\in \Sigma$, which says that $x\star x=0$ and $x=0$ as $(A,\star)$ is semi-prime.
\end{proof}
\end{corollary}

\section{Tensor Products}

We actually prove rather more than that there exists a tensor product of Archimedean $f$-algebras $A$ and $B$, namely that the Fremlin positive projective tensor product of $A$ and $B$ itself, $A\ftp B$, can be given an $f$-algebra multiplication that does what is required. We refer the interested reader to Fremlin's original paper \cite{F} for the basics of this construction. The basis for our construction is a simple description of the Fremlin tensor product using representations in $S(\Sigma)$.

\begin{proposition}
Let $E$ and $F$ be vector sublattices of $S(\Sigma)$ and $S(\Omega)$ respectively. Then the Fremlin tensor product of $E$ and $F$, $E\ftp F$ is vector lattice isomorphic to the vector sublattice of $S(\Sigma\times\Omega)$ generated by the functions of the form $(\sigma,\omega)\mapsto x(\sigma) y(\omega)$ for $x\in E$ and $y\in F$.
\begin{proof}
The bilinear map $\psi:E\times F\to S(\Sigma\times \Omega)$ defined by the property $\psi(x,y)(\sigma,\omega)=x(\sigma) y(\omega)$ certainly has the property that if $0\lneqq x\in E$ and $0\lneqq y\in F$ then $0\lneqq \psi(x,y)$. The claim now follows from Corollary 4.4 of \cite{F}.
\end{proof}
\end{proposition}

\begin{remark}
Note that any Archimedean vector lattices may be represented in this way. Note also that this argument does most emphatically \emph{not} give a simple construction of the Fremlin tensor product as it depends critically on Fremlin's work in \cite{F}.
\end{remark}

If $(A,\star)$ and $(B,\bullet)$ are Archimedean $f$-algebras then the algebraic tensor product, $A\otimes B$, has a natural algebraic multiplication $\times$ characterized by
\[(a\otimes b)\times (a'\otimes b')=(a\star a')\otimes (b\bullet b').\]
Unfortunately, in general $A\otimes B$ is much smaller than $A\ftp B$.

\begin{theorem}\label{ftensor}
If $A$ and $B$ are Archimedean $f$-algebras then the multiplication $\times$ on $A\otimes B$ has a unique extension to an $f$-algebra multiplication on the  Fremlin positive tensor product $A\ftp B$.
\begin{proof}
Using Proposition \ref{frep2} and any suitable representation theorem for Archimedean vector lattices, we may identify $A$ and $B$ with admissible vector sublattices of $S(\Sigma)$ and $S(\Omega)$ respectively, with there being $v\in S(\Sigma)_+$ and $w\in S(\Omega)_+$ such that $(a\star a')(\sigma)=v(\sigma) a(\sigma) a'(\sigma)$ and $(b\bullet b')(\omega)=w(\omega) b(\omega) b'(\omega)$ whenever the pointwise products are defined. On the algebraic tensor product, we have
\begin{align*}\big((a\otimes b)\star (a'\otimes b')\big)(\sigma,\omega)&=\big((a\star a')\otimes (b\bullet b')\big)(\sigma,\omega)\\
&=(a \star a')(\sigma)(b\bullet b')(\omega)\\
&=\big(v(\sigma) a(\sigma)a'(\sigma)\big)\big(w(\omega) b(\omega) b'(\omega)\big)\\
&=\big(v(\sigma)w(\omega)\big) (a\otimes b)(\sigma,\omega) (a'\otimes b')(\sigma,\omega).
\end{align*}
Let us  define $u(\sigma,\omega)$ to be (the equivalence class of) the function $v(\sigma)w(\omega)$ on the Cartesian product $U\times V$ of dense open subsets $U\subset \Sigma$ and $V\subset \Omega)$ on which $v$ and $w$, respectively, are defined. As $U\times V$ is certainly a dense open subset of $\Sigma\times \Omega$, it is clear that $u\in S(\Sigma\times\Omega)_+$. It follows from Corollary \ref{frep2} that
\[f\times g(\sigma,\omega)=u(\sigma,\omega) f(\sigma,\omega) g(\sigma,\omega)\]
defines an $f$-algebra multiplication on $S(\Sigma, \Omega)$, which extends the original multiplication $\times$ on $A\times B$. By  Theorem 2.1 (B) of \cite{HaJ}, the sublattice of $S(\Sigma,\Omega)$ generated by $A\otimes B$, which is precisely $A\ftp B$, is an $f$-algebra under the multiplication $\times$.

To see the uniqueness, we recall that multiplication by a positive element in an $f$-algebra is a lattice homomorphism. Using the fact that any element of $A\ftp B$ is a finite supremum of a finite infimum of elements of $A\otimes B$, we see that any $f$-algebra extension of the multiplication $\times$ on $A\otimes B$ must coincide with our product when $a\in (A\otimes B)_+$ and $a'\in (A\ftp B)$. Repeat this argument to see that we have equality when $a\in A\ftp B$ and $a'\in (A\ftp B)_+$ and the final step is straightforward.
\end{proof}
\end{theorem}

In view of the canonical nature of this construction, we will henceforth refer to it as \emph{the} $f$-algebra tensor product of $A$ and $B$.
As is well known, nice order theoretic properties of $A$ and $B$ are seldom inherited by $A\ftp B$. It turns out, however, that the situation is rather better for algebraic properties.

\begin{theorem}
The $f$-algebra tensor product $(A\ftp B,\times)$ of non-zero $f$-algebras $(A,\star)$ and $(B,\bullet)$ is semi-prime if and only if both $(A,\star)$ and $(B,\bullet)$ are semi-prime.
\begin{proof}
Recall that an $f$-algebra is semi-prime if and only if $x^2=0\Rightarrow x=0$.  If, for example, $A$ fails to be semi-prime there is $0\ne a\in A$ with $a\star a=0$. Taking any non-zero $b\in B$, $0\ne a\otimes b\in A\ftp B$ but $(a\otimes b)\times (a\otimes b)=(a\star a)\otimes (b\bullet b)=0\otimes (b\bullet b)$ is zero so that $A\ftp B$ is not semi-prime.

If $A$ and $B$ are semi-prime represent $(A,\star)$ in $S(\Sigma)$ and $(B,\bullet)$ in $S(\Delta)$ with the multiplication weights $v\in S(\Sigma)_+$ and $w\in S(\Delta)_+$ being strictly positive, using Corollary \ref{semiprimerep2}. The construction of the multiplication in $A\ftp B$ in Theorem \ref{ftensor} involved the weight $(\sigma,\delta)\mapsto v(\sigma) w(\delta)$ which is strictly positive, so the remark after Corollary \ref{frep2} tells us that $(A\ftp B,\times)$ is semi-prime.
\end{proof}
\end{theorem}

Recall that if $X$ is a vector subspace of a vector lattice then $X^{\vee}$ is the set of all finite suprema from $X$, $X^{\wedge}$ is the set of all finite infima from $X$ and that the vector sublattice generated by $X$ is precisely $X^{\vee\wedge}=X^{\wedge\vee}$.

\begin{theorem}
The $f$-algebra tensor product of $f$-algebras $A$ and $B$ has a multiplicative identity if and only if both $A$ and $B$ have multiplicative identities.
\begin{proof}
If $A$ and $B$ have multiplicative identities $e_A$ and $e_B$ respectively then $e_A\otimes e_B$ is a multiplicative identity on $A\otimes B$ and it is clear, given that in an $f$-algebra multiplication by positive elements is a lattice homomorphism, that it remains a multiplicative identity in $A\ftp B$.

Now suppose that $A\ftp B$ has a multiplicative identity $e$, and is  hence semi-prime. Then $A$ and $B$ are semi-prime and we may represent $(A,\star)$ in $S(\Sigma)$, with strictly positive weight $v$, and $(B,\bullet)$ in $S(\Delta)$, with strictly positive weight $w$ and identify $(A\ftp B,\times)$ with the vector sublattice of $S(\Sigma\times \Delta)$ generated by $A\otimes B$ with the multiplication being given by the weight $u(\sigma,\delta)=v(\sigma)w(\delta)$ (on the product of dense open sets on which $v$ and $w$ are defined.)

As $e$ is a multiplicative identity we certainly have $e\times(a\otimes b)=a\otimes b$ for all $a \in A$ and $b \in B$. Let $e$ be defined (at least) on the dense open set $U\subset \Sigma\times \Delta$. For any $(\sigma,\delta)\in U$ we may choose $a\in A $ and $b\in B$ with $a(\sigma)\ne 0\ne b(\delta)$. As
\[v(\sigma)w(\delta) e(\sigma,\delta)a(\sigma)b(\delta)=(e\times (a\otimes b)(\sigma,\delta)=(a\otimes b)(\sigma,\delta)=a(\sigma)b(\delta)\]
we see that $e(\sigma,\delta)=\big(v(\sigma)w(\delta)\big)^{-1}$ on $U$. As $\big(v(\sigma)w(\delta)\big)^{-1}$ is actually defined on the whole of $\Sigma\times \Delta$ (remembering that both weights are strictly positive) we may take $U$ to be the whole of $\Sigma\times \Delta$. It will be  clear from the representations of the multiplications that $e_A(\sigma)=v(\sigma)^{-1}$ and $e_B(\delta)=w(\delta)^{-1}$ are multiplicative identities for $A$ and $B$ respectively, \emph{once we show that they do lie in the appropriate space.}

We do this by showing that any element $h\in A\ftp B$ is defined on a set of the form $U\times V$, where $U$ (resp. $V$) is a dense open subset of $\Sigma$ (resp. $\Delta$) and that, for example, if $\delta\in V$ then $\sigma\mapsto h(\sigma,\delta)$ is an element of $A$ (defined at least on $U$.) The claim that $e_A\in A$ and that $e_B\in B$ will then be obvious. Suppose first that $h\in A\otimes B$ then we may write $h=\sum_{k=1}^n f_k\otimes g_k$ with each $f_k$ defined on a dense open $U_k\subset \Sigma$ and $g_k$ defined on a dense open $V_k\subset \Delta$. Then $h$ is certainly defined on $U\times V$ where $U=\bigcap_{k=1}^n U_k$ and $V=\bigcap_{k=1}^n V_k$ are dense open subsets of $\Sigma$ and $\Delta$ respectively. For any $\delta\in V$, $h(\sigma,\delta)=\sum_{k=1}^n g_k(\delta)f_k(\sigma)$ is just a finite linear combination of elements $f_k\in A$, so lies in $A$. Given this, a similar argument shows that if $h\in (A\otimes B)^{\wedge}$ then $h$ has the appropriate domain and that $\sigma\mapsto h(\sigma,\delta)$ is a finite infimum of elements from $A$ so lies in $A$. Finally a third step shows the same for elements of $(A\otimes B)^{\wedge\vee}=A\ftp B$ and the proof is complete.
\end{proof}
\end{theorem}

It is clear that there is no other sensible way to define an $f$-algebra structure on $A\ftp B$, but we may still ask what nice categorical properties does it have?
In a purely algebraic setting the normal universal property for a tensor product of two algebras $A$ and $B$ is that if $T_A:A\to C$ and $T_B:B\to C$ are algebra homomorphisms with commuting ranges then they induce an algebra homomorphism $S:A\otimes B\to C$ with $S(a\otimes b)=T_A(a) T_B(b)$. Often there is some assumption about identity elements as well. In our setting the commutativity of the ranges is automatic, but we should take the order structure into account in some way. Order boundedness of the algebra homomorphisms is about the least that we could ask. Bearing in mind Theorem 5.1 of \cite{HdP}, which tells us that both the domain and range are semi-prime and the domain is uniformly complete then every $f$-algebra homomorphism is a lattice homomorphism, this is not such a strong restriction. We do not, thankfully, need to assume anything about identities so that we do not need to open the can of worms of unital embeddings of $f$-algebras.

Much of the work of proving the required universal property is contained in the following specialized, but quite general, lemma.  We also need the fact that in any $f$-algebra $(A,\star)$ if $x\ge 0$ then $x\star (y\vee z)=(x\star y)\vee (x\star z)$ and $x\star (y\wedge z)=(x\star y)\wedge (x\star z)$ for any $ y,z\in A$.

\begin{lemma}\label{multext}
Let $(A,\star)$ and $(B,\bullet)$ be $f$-algebras, $X$ a sub-algebra of $A$ such that $X=X\cap A_+-X\cap A_+$, $T:A\to B$ a lattice homomorphism such that $T_{|X}$ is an algebra homomorphism, then $T_{|X^{\wedge\vee}}$ is also an algebra homomorphism.
\begin{proof}
We start by proving that if $Y$ and $Z$ are positively generated subspaces of $A$ with
\[T(y\star z)=Ty\bullet Tz\eqno{(\dag)}\]
for all $y\in Y$ and all $z\in Z$ then also ($\dag$) holds for all $y\in Y$ and $z\in Z^\wedge$.
If $z_k\in Z$ for $1\le k\le n$ and $y\in Y\cap A_+$ then
\begin{align*}
Ty\bullet T(\wedge_{k=1}^n z_k)&=Ty\bullet \wedge_{k=1}^n Tz_k\\
&=\wedge_{k=1}^n Ty\bullet Tz_k\\
&=\wedge_{k=1}^n T(y\star z_k)\\
&=T(\wedge_{k=1}^n y\star z_k)\\
&=T(y\star \wedge_{k=1}^n z_k)
\end{align*}
and then linearity gives the same conclusion for any $y\in Y$. A similar argument proves the corresponding result for $z\in Z^\vee$, Now, starting with $Y=Z=X$ we see that ($\dag$) holds in succession for $y\in X, z\in X^\wedge$; $y\in X^\wedge, z\in X^\wedge$; $y\in X^\wedge, z\in X^{\wedge\vee}$ and finally for $y\in X^{\wedge\vee}, z\in X^{\wedge\vee}$.
\end{proof}
\end{lemma}

\begin{theorem}If $(A,\star), (B,\circledast)$ and $(C,\bullet)$ are $f$-algebras, $T_A:A\to C$ and $T_B:B\to C$ are order bounded algebra homomorphisms then there is a unique positive algebra and lattice  homomorphism $S:A\ftp B\to C$ with $S(a\otimes b)=T_A(a)\bullet T_B(b)$ for all $a\in A$ and $b\in B$.
\begin{proof}

Let $N=\{c\in C:c\bullet c=0\}$ which is a band and algebra ideal in $C$. The quotient $C/N$ is naturally a semi-prime $f$-algebra. Let $Q:C\to C/N$ be the quotient map. We know (see  Proposition 10.2 of \cite{dP}) that all products $c\bullet c'$ lie in $N^d$, the band complementary to $N$. Both $Q\circ T_A$ and $Q\circ T_B$ are order bounded algebra homomorphisms into the semi-prime $f$-algebra $C/N$ so are actually lattice homomorphisms by Theorem 4.5 of \cite{T}.
 The map $(a,b)\mapsto T_A(a)\bullet T_B(b):A\times B\to N^d\subset C$ is now easily seen to be a lattice bimorphism, as $Q:N^d\to C/N$ is a one-to-one $f$-algebra homomorphism and $Q\circ \big( T_A(a)\bullet T_B(b)\big)=\big(Q\circ T_A(a)\big)\bullet \big(Q\circ T_B(b)\big)$. Now by Theorem 4.2 of \cite{F} there is a a linear lattice homomorphism $S:A\ftp B\to N^d\subset C$ with $S(a\otimes b)=T_A(a)\bullet T_B(b)$. Routine algebra shows that, on the algebraic tensor product $A\otimes B$,  $S$ is an algebra homomorphism. As $A\otimes B$ is a positively generated subspace of $A\ftp B$, Lemma \ref{multext} shows that $T$ is actually an algebra homomorphism on $(A\otimes B)^{\wedge\vee}=A\ftp B$. The uniqueness follows from $S$ being a lattice homomorphism that is specified on the lattice generating subset $A\otimes B$ or $A\ftp B$.
\end{proof}
\end{theorem}


\begin{bibdiv}
\begin{biblist}[\resetbiblist{99}]

\bib{AAJ}{article}{
   author={Azouzi, Y.},
   author={Ben Amor, M. A.},
   author={Jaber, J.},
   title={The tensor product of function algebras},
  journal={arXiv:1512.00703v1},
}
\bib{BP}{article}{
   author={Birkhoff, Garrett},
   author={Pierce, R. S.},
   title={Lattice-ordered rings},
   journal={An. Acad. Brasil. Ci.},
   volume={28},
   date={1956},
   pages={41--69},
   issn={0001-3765},
   review={\MR{0080099 (18,191d)}},
}
\bib{BvR}{article}{
   author={Buskes, G.},
   author={van Rooij, A.},
   title={Almost $f$-algebras: structure and the Dedekind completion},
   note={Positivity and its applications (Ankara, 1998)},
   journal={Positivity},
   volume={4},
   date={2000},
   number={3},
   pages={233--243},
   issn={1385-1292},
   review={\MR{1797126 (2001j:46062)}},
   doi={10.1023/A:1009874426887},
}

\bib{F}{article}{
   author={Fremlin, D. H.},
   title={Tensor products of Archimedean vector lattices},
   journal={Amer. J. Math.},
   volume={94},
   date={1972},
   pages={777--798},
   issn={0002-9327},
   review={\MR{0312203 (47 \#765)}},
}

\bib{GL}{article}{
   author={Grobler, J. J.},
   author={Labuschagne, C. C. A.},
   title={An $f$-algebra approach to the Riesz tensor product of Archimedean
   Riesz spaces},
   journal={Quaestiones Math.},
   volume={12},
   date={1989},
   number={4},
   pages={425--438},
   issn={0379-9468},
   review={\MR{1021941 (91d:46005)}},
}
\bib{HaJ}{article}{
   author={Hager, A. W.},
   author={Johnson, D. G.},
   title={Some comments and examples on generation of (hyper-)archimedean
   $\ell$-groups and $f$-rings},
   language={English, with English and French summaries},
   journal={Ann. Fac. Sci. Toulouse Math. (6)},
   volume={19},
   date={2010},
   number={Fascicule Special},
   pages={75--100},
   issn={0240-2963},
   review={\MR{2675722 (2011i:06032)}},
}
\bib{HeJ}{article}{
   author={Henriksen, M.},
   author={Johnson, D. G.},
   title={On the structure of a class of archimedean lattice-ordered
   algebras. },
   journal={Fund. Math.},
   volume={50},
   date={1961/1962},
   pages={73--94},
   issn={0016-2736},
   review={\MR{0133698 (24 \#A3524)}},
}


\bib{HdP}{article}{
   author={Huijsmans, C. B.},
   author={de Pagter, B.},
   title={Subalgebras and Riesz subspaces of an $f$-algebra},
   journal={Proc. London Math. Soc. (3)},
   volume={48},
   date={1984},
   number={1},
   pages={161--174},
   issn={0024-6115},
   review={\MR{721777 (85f:46015)}},
   doi={10.1112/plms/s3-48.1.161},
}
		



\bib{J}{article}{
   author={Johnson, D. G.},
   title={On a representation theory for a class of Archimedean
   lattice-ordered rings},
   journal={Proc. London Math. Soc. (3)},
   volume={12},
   date={1962},
   pages={207--225},
   issn={0024-6115},
   review={\MR{0141685 (25 \#5082)}},
}
\bib{MN}{book}{
   author={Meyer-Nieberg, Peter},
   title={Banach lattices},
   series={Universitext},
   publisher={Springer-Verlag, Berlin},
   date={1991},
   pages={xvi+395},
   isbn={3-540-54201-9},
   review={\MR{1128093 (93f:46025)}},
   doi={10.1007/978-3-642-76724-1},
}

\bib{dP}{thesis}{
   author={de Pagter, B.},
   title={$f$-algebras and orthomorphisms},
      type={Ph.D. thesis},
      organization={University of Leiden},
   date={1981},
   pages={145},

}

\bib{T}{article}{
   author={Triki, Abdelmajid},
   title={On algebra homomorphisms in complex almost $f$-algebras},
   journal={Comment. Math. Univ. Carolin.},
   volume={43},
   date={2002},
   number={1},
   pages={23--31},
   issn={0010-2628},
   review={\MR{1903304}},
}

\bib{W}{article}{
   author={Wickstead, A. W.},
   title={Representation and duality of multiplication operators on
   Archimedean Riesz spaces},
   journal={Compositio Math.},
   volume={35},
   date={1977},
   number={3},
   pages={225--238},
   issn={0010-437X},
   review={\MR{0454728 (56 \#12976)}},
}

\end{biblist}
\end{bibdiv}
\end{document}